\documentclass[11pt,bezier]{article}
\usepackage{amsmath,amssymb,amsfonts,euscript}

\usepackage{graphicx}

\usepackage[ansinew]{inputenc}
\usepackage{graphicx}
\usepackage{color}
\usepackage{mathrsfs}

\usepackage[colorlinks]{hyperref}
\usepackage[active,new,noold,marker]{xrcs}
\catcode`\=13
\def{$\bowtie$}

\usepackage{eurosym}
\usepackage{lscape} %landcape pages support

\newtheorem{prethm}{{\bf Theorem}}

\newenvironment{thm}{\begin{prethm}{\hspace{-0.5
               em}{\bf}}}{\end{prethm}}

\newtheorem{prepro}{{\bf Theorem}}

\newtheorem{preprop}{{\bf Proposition}}

\newtheorem{precor}{{\bf Corollary}}

\newtheorem{preconj}{{\bf Conjecture}}

\newtheorem{preconj1}{{\bf Conjecture}}

\newenvironment{conj1}{\begin{preconj1}{\hspace{-0.5
               em}{\bf}}}{\end{preconj1}}

\newtheorem{predefi}{{\bf Definition}}

\newtheorem{preremark}{{\bf Remark}}

\newtheorem{preexample}{{\bf Fact}}

\newenvironment{example}{\begin{preexample}\rm{\hspace{-0.5
               em}{\bf}}}{\end{preexample}}

\newtheorem{prelem}{{\bf Lemma}}

\newtheorem{prelam}{{\bf Lemma}}

\newtheorem{preprob}{{\bf Problem}}

\newtheorem{preproof}{{\bf Proof}}

\newtheorem{preali}{{\bf Proof of Theorem 1.}}

\newenvironment{ali}[1]{\begin{preali}{\rm
               #1}\hfill{$\Box$}}{\end{preali}}

\newtheorem{prealii}{{\bf Proof of Theorem 2.}}

\newenvironment{alii}[1]{\begin{prealii}{\rm
               #1}\hfill{$\Box$}}{\end{prealii}}

\newtheorem{prealiii}{{\bf Proof of Theorem 3.}}

\newenvironment{aliii}[1]{\begin{prealiii}{\rm
               #1}\hfill{$\Box$}}{\end{prealiii}}

\newtheorem{prealiiii}{{\bf Proof of Theorem 4.}}

\newenvironment{aliiii}[1]{\begin{prealiiii}{\rm
               #1}\hfill{$\Box$}}{\end{prealiiii}}

\newtheorem{prealiiiii}{{\bf Proof of Theorem 5.}}

\newenvironment{aliiiii}[1]{\begin{prealiiiii}{\rm
               #1}\hfill{$\Box$}}{\end{prealiiiii}}

\title{ The inapproximability for the  (0,1)-additive number}

\author{{\normalsize
{Arash Ahadi${}^{\mathsf{a}}$},\,
 {Ali Dehghan${}^{\mathsf{b}}$},\,
}\vspace{3mm}
\\{\footnotesize{${}^{\mathsf{a}}$\it Department of
Mathematical Sciences, Sharif University of Technology, Tehran,
Iran}}  {\footnotesize{}}\\{\footnotesize{${}^{\mathsf{b}}$\it
Systems and Computer Engineering Department,
Carleton University, Ottawa, Ontario, Canada}}
\thanks{{\it E-mail addresses}:  $\mathsf{arash\_ahadi@mehr.sharif.edu}$ (Arash Ahadi),  $\mathsf{ali\_dehghan16@aut.ac.ir}$ (Ali Dehghan). } }

\date{}
\begin{document}
\maketitle

\begin{abstract}
{\small \noindent
An
 {\it additive labeling} of a graph $G$  is a function $ \ell :V(G) \rightarrow\mathbb{N}$, such that for every two adjacent vertices $ v $ and $
u$ of $ G $, $  \sum_{w \sim v}\ell(w)\neq  \sum_{w \sim u}\ell(w) $ ($ x \sim y $
means that $ x $ is joined to $y$).
The {\it additive number} of $ G $, denoted by $\eta(G)$, is the minimum number $k $ such that $ G $ has a additive labeling $ \ell :V(G) \rightarrow \mathbb{N}_k$.
The
{\it additive choosability} of a graph $G$, denoted by $\eta_{\ell}(G) $,
is the
smallest number $k$ such that $G$ has an additive labeling for any assignment
of lists of size $k$ to the vertices of $G$, such that the label of each vertex
belongs to its own list.

Seamone (2012) \cite{a80} conjectured that for every graph $G$, $\eta(G)= \eta_{\ell}(G)$.
We give a negative answer to this
conjecture  and we  show that  for
every $k$  there is a graph $G$ such that $ \eta_{\ell}(G)- \eta(G) \geq k$.

A
{\it $(0,1)$-additive labeling} of a graph $G$ is a function $ \ell :V(G) \rightarrow\{0,1\}$, such that for every two adjacent vertices $ v $ and $
u$ of $ G $, $ \sum_{w \sim v}\ell(w)\neq \sum_{w \sim u}\ell(w) $.
 A graph may lack any $(0,1)$-additive labeling. We show that it is $ \mathbf{NP} $-complete to decide whether a $(0,1)$-additive labeling exists for some families of graphs such as perfect graphs and planar triangle-free graphs.
For a graph $G$ with some $(0,1)$-additive labelings, the $(0,1)$-additive number of $G$ is defined as $ \sigma_{1} (G) = \min_{\ell \in \Gamma}\sum_{v\in V(G)}\ell(v)  $ where $\Gamma$ is the set of $(0,1)$-additive labelings of $G$.
We prove that given a planar graph that admits a $(0,1)$-additive labeling, for all $
\varepsilon >0 $, approximating the $(0,1)$-additive number
within $ n^{1-\varepsilon} $ is $ \mathbf{NP} $-hard.
}

\begin{flushleft}


\end{flushleft}

\end{abstract}

\section{Introduction}
\label{}
Throughout the paper we denote $\{1,2,\ldots, k\}$ by $\mathbb{N}_k$.
An
{\it additive labeling} of a graph $G$, which was introduced by Czerwi\'nski et al. \cite{MR2552893}, is a function $ \ell :V(G) \rightarrow\mathbb{N}$, such that for every two adjacent vertices $ v $ and $
u$ of  $ G $, $  \sum_{w \sim v}\ell(w)\neq  \sum_{w \sim u}\ell(w) $ ($ x \sim y $
means that $ x $ is joined to $y$).
The {\it additive  number} of $ G $, denoted by $\eta(G)$, is the minimum number $k $ such that $ G $ has a additive labeling $ \ell :V(G) \rightarrow \mathbb{N}_k $.
Initially, additive labeling
is called a lucky labeling of $G$.
The following important conjecture was proposed by Czerwi\'nski et al. \cite{MR2552893}.

\begin{conj1}
 {\bf[ Additive Coloring Conjecture \rm \cite{MR2552893}]\bf} For every
graph $G$, $\eta(G)\leq \chi(G)$.
\end{conj1}

Czerwi\'nski et al. also, considered the list version of above problem \cite{MR2552893}.
The
{\it additive choosability} of a graph $G$, denoted by $\eta_{\ell}(G) $, is
 the smallest number $k $ such that $G$ has an additive labeling from any
assignment of lists of size $k$ to the vertices of $G$.
Idem above, about list-coloring proved that if $T$ is a tree, then $\eta_{\ell}(T)\leq 2 $, and if $G$ is a bipartite planar graph,
then $\eta_{\ell}(G)\leq 3 $ (for more information about the recent results see \cite{brandt2015lucky}). Seamone in his
Ph.D dissertation posed the following conjecture about the relationship between additive  number and additive choosability \cite{a80}.

\begin{conj1}
 {\bf [Additive List Coloring Conjecture \rm \cite{a80}]\bf} For every
graph $G$, $\eta(G)= \eta_{\ell}(G)$.
\end{conj1}

For a given connected graph $G$ with at least two vertices, if no two adjacent
vertices have a same degree, then $\eta(G)=1$ and $\eta_{\ell}(G)>1$.
We show that not only there exists a counterexample for the above equality but also the difference between
$\eta(G)$ and $ \eta_{\ell}(G)$ can be arbitrary large.

\begin{thm}
 For every $k$ there is a graph $G$ such that $  \eta(G) \leq k \leq \eta_{\ell}(G)/2$.
\end{thm}

 Chartrand et al. introduced another version of additive labeling and called it sigma coloring
\cite{MR2729020}. For a graph $G$, let $ c :V(G) \rightarrow \mathbb{N}$ be a vertex labeling of
$G$.  If for every two adjacent vertices $ v $ and $
u$ of $ G $, $  \sum_{w \sim v}c(w)\neq  \sum_{w \sim u}c(w) $, then $c$ is called a {\it sigma coloring} of $G$. The minimum number of labels
required in a sigma coloring is called the {\it sigma chromatic number} of $G$
and is denoted by $\sigma(G)$. Chartrand et al. proved that, for every
graph $G$, $\sigma(G) \leq \chi(G)$ \cite{MR2729020}.
Note that the only difference between additive labeling and sigma coloring is the
objective function, but the feasible labelings are the same.

Additive labeling and sigma coloring have been studied extensively by several
authors, for instance see \cite{ ahadi, akbari2, additive,   aaa, MR2729020,  MR2552893, ali,  MR2729390, marenco2015topological}.
It is proved, in \cite{ahadi} that it is $ \mathbf{NP} $-complete to determine whether a given graph $G$ has $ \eta(G)=k$
for any $k\geq 2$. Also,
it was shown that, it is $ \mathbf{NP} $-complete to decide for a given planar 3-colorable graph $G$, whether $ \eta(G)=2$ \cite{ahadi}. Furthermore, it was proved that, it is $ \mathbf{NP} $-complete to decide  for a given 3-regular graph $G$, whether $ \eta(G)=2$ \cite{ali}.

The edge version of additive labeling was introduced by
Karo\'nski, \L{}uczak and Thomason \cite{MR2047539}. They introduced an edge-labeling which
is additive vertex-coloring that means for every edge
$uv$, the sum of labels of the edges incident to $u$ is different
from the sum of labels of the edges incident to $v$ \cite{MR2047539}.
It is conjectured that three integer labels
$\mathbb{N}_3$ are sufficient for every connected graph, except
$K_2$ \cite{MR2047539}. Currently the best bound is $5$ \cite{MR2595676}.
This labeling has been studied extensively by several
authors, for instance see \cite{MR2145514,   MR2404230, z2, dehghan2013algorithmic, z3, D2}.

A {\it clique} in a  graph $G = (V, E)$  is a subset of its vertices such that every two vertices in the subset are connected by an edge.
The {\it clique number} $\omega(G)$ of a graph $G$ is the number of vertices in a maximum clique in $G$.
There is no direct relationship between the additive number  and the clique number of graphs. For any natural number $\omega$ there exists a graph $G $, such that $ \omega(G )=\omega $ and $ \eta(G )=1  $. To see this for given number $ \omega $, consider a graph $G $ with the set of vertices
$V(G )= \lbrace v_{i} \vert  i \in \mathbb{N}_{\omega} \rbrace \cup \lbrace u_{i,j} \vert    i,j\in \mathbb{N}_{\omega}, j < i   \rbrace $
and the set of edges
$ E(G )=\lbrace v_{i}v_{j}\vert i\neq j \rbrace \cup \lbrace v_{i}u_{i,j}\vert  i,j\in \mathbb{N}_{\omega}, j < i   \rbrace $.

\begin{thm} We have the following:\\
{\em (i)} For every graph $G$, $ \eta (G)\geq \frac{w}{n-w+1}$.\\
{\em (ii)} If $ G $ is a regular graph and $ \omega > \frac{n+4}{3} $, then $ \eta (G) \geq  3 $.
\end{thm}

A
{\it $(0,1)$-additive labeling} of a graph $G$ is a function $ \ell :V(G) \rightarrow\{0,1\}$, such that for every two adjacent vertices $ v $ and $
u$ of $ G $, $  \sum_{w \sim v}\ell(w)\neq  \sum_{w \sim u}\ell(w) $.
A graph may lack any $(0,1)$-additive labeling.  It was proved that, it is $ \mathbf{NP} $-complete to decide  for a given 3-regular graph $G$, whether $ \eta(G)=2$ \cite{ali}. So, it is $ \mathbf{NP} $-complete to decide whether a $(0,1)$-additive labeling exists for a given 3-regular graph $G$. In this paper, we study the computational complexity of $(0,1)$-additive labeling for perfect graphs and planar graphs.

A graph $G$ is called {\it perfect} if $\omega(H) = \chi(H)$ for every
induced subgraph $H$ of $G$.
Here, we show that it is $ \mathbf{NP} $-complete to decide whether a $(0,1)$-additive labeling exists for perfect graphs.

\begin{thm}
The following problem is $ \mathbf{NP} $-complete:
 Given a perfect graph $G$, does $G$ have any $(0,1)$-additive labeling?
\end{thm}

Next, we show that it is $ \mathbf{NP} $-complete to decide whether a $(0,1)$-additive labeling exists for  planar triangle-free graphs.

\begin{thm}
It is $ \mathbf{NP} $-complete to determine whether a given  planar triangle-free graph $G$ has a $(0,1)$-additive labeling.
\end{thm}

For a graph $G$ with some $(0,1)$-additive labelings, the $(0,1)$-additive number of $G$ is defined as $ \sigma_{1} (G) = \min_{\ell \in \Gamma}  \sum_{v\in V(G)}\ell(v)  $ where $\Gamma$ is the set of $(0,1)$-additive labelings of $G$. For a given graph $G$ with  a $(0,1)$-additive labeling $\ell$ the function $f(v)=1+\sum_{w \sim v}\ell(w)$ is a proper vertex coloring, so we have the following trivial lower bound for $ \sigma_{1} (G)$.

\begin{center}
$\chi(G) - 1 \leq \sigma_{1} (G)$.
\end{center}

We prove that given a planar graph that admits a $(0,1)$-additive labeling, for all $
\varepsilon >0 $, approximating the $(0,1)$-additive number
within $ n^{1-\varepsilon} $ is $ \mathbf{NP} $-hard.

\begin{thm}
If $ \mathbf{P} \neq \mathbf{NP}$, then for any constant $\varepsilon > 0$, there
is no polynomial-time $ n^{1-\varepsilon} $-approximation algorithm for finding $ \sigma_{1} (G) $ for a given planar graph with at least one $(0,1)$-additive labeling.
\end{thm}

 For $v \in V (G)$ we denote by $N(v)$ the set of
neighbors of $v$ in $G$. Also, for every $v \in V (G)$,  the degree of $v$ is  denoted by $d(v)$.
 We follow \cite{MR1567289, MR1367739} for terminology and
notation not defined here, and we consider finite undirected
simple graphs $G=(V,E)$.

\section{Counterexample}
 \label{T1}

\begin{ali}{

 For
every $k$ we construct a graph $G$ such that $ \eta_{\ell}(G)- \eta(G) \geq k$.
For every $\alpha$, $ \alpha \in \mathbb{N}_{2k-1} $ consider a copy of complete graph $K^{(\alpha)}_{2k}$, with the vertices $\{ x_{\beta}^{\alpha} :   \beta \in \mathbb{N}_k  \}  \cup \{y_{\beta}^{\alpha} :   \beta \in \mathbb{N}_k  \}$. Next, consider an isolated vertex $t$ and join every vertex $y_{\beta  }^{\alpha}$ to $t$,
Call the resulting graph $G$.
First, note that in every additive labeling $\ell$ of $G$, for
every $(i,j)$, where $ i < j $ and $i,j\in \mathbb{N}_k$ we have $\sum_{z\in N(x_{i}^{1})}\ell(z)\neq \sum_{z\in N(x_{j}^{1})}\ell(z) $, thus
$\ell(x_{i}^{1})\neq \ell(x_{j}^{1})$ (because all the neighbors of $x_{i}^{1}$ and $x_{j}^{1}$
are common except $x_{i}^{1}$ as a neighbor of $x_{j}^{1}$, and vice
versa). Therefore $\ell(x_{1}^{1}), \ell(x_{2}^{1}),\ldots, \ell(x_{k}^{1}) $ are $k$ distinct numbers, that means $\eta(G) \geq k$.
Define (for every $\alpha$ and $\beta$):

 $\ell: V(G)\rightarrow \mathbb{N}_{k}$,

 $\ell(x_{\beta}^{\alpha})=\ell(y_{\beta}^{\alpha})=\beta$,

 $\ell(t)=k$.

 It is easy to see that $\ell$ is an additive labeling for $G$. Next, we show that $\eta_{\ell}(G) > 2k-1$.
 Consider the following lists for the vertices of $G$ (for every $\alpha$ and $\beta$).

 $L(x_{\beta}^{\alpha})=\mathbb{N}_{2k-1}$,

$L(y_{\beta}^{\alpha})=\{i+\alpha : i \in \mathbb{N}_{2k-1} \}$,

 $L(t)=\mathbb{N}_{2k-1}$.

 To the contrary suppose that $\eta_{\ell}(G) \leq 2k-1$ and let $\ell$ be an additive labeling from the above lists. Suppose that $\ell(t)=r$. Consider the complete graph $K^{(r)}_{2k}$, for every   $\beta$ we have:

 $L(x_{\beta}^{r})=\mathbb{N}_{2k-1}$,   

$L(y_{\beta}^{r})=\{i+r : i \in \mathbb{N}_{2k-1} \}$.

 Now, consider the following partition for  $\mathbb{N}_{2k-1}\cup \{i+r : i \in \mathbb{N}_{2k-1} \}$,

 \begin{center}
$\{1+r,1\},\{2+r,2\}, \ldots, \{2k-1+r,2k-1\}$.
 \end{center}

 By Pigeonhole Principle, there are indices $i$, $n$ and $m$ such that $\ell(x_{m}^{r}),\ell(y_{n}^{r})\in \{i+r,i\}$, so $\ell(x_{m}^{r})=i$ and $\ell(y_{n}^{r})=i+r$. Therefore, $\sum_{z\in N(x_{m}^{r})}\ell(z)= \sum_{z\in N(y_{n}^{r})}\ell(z) $. This is a contradiction, so $\eta_{\ell}(G) \geq 2k$.

}\end{ali}

\section{Lower bounds}

\begin{alii}{
{\bf (i)} Let $ \ell:V(G)\rightarrow \mathbb{N}_k $ be an additive labeling of $ G $ and suppose that $ T=\lbrace v_i | i \in \mathbb{N}_{\omega} \rbrace$ is a maximum clique in $ G $. For each vertex $ v\in T $, define the function $Y_{v}$.

\begin{center}
 $ Y_{v}\stackrel {\rm def}{=}{ \displaystyle\sum_{x\in V(G)\setminus T \atop x\sim v} l(x) -l(v)} $.
\end{center}

For every two adjacent vertices $ v$ and $u $ in $ T $, we have:

\begin{center}

$\displaystyle\sum_{ x\sim v} l(x)\neq \sum_{ x\sim u} l(x)$,

$  \displaystyle\sum_{x\notin T \atop x\sim v} l(x)+ \sum_{x\in T \atop x\neq v} l(x) \neq   \sum_{x\notin T \atop x\sim u} l(x)+ \sum_{x\in T \atop x\neq u} l(x)$,

$  \displaystyle\sum_{x\notin T \atop x\sim v} l(x)+ l(u) \neq  \sum_{x\notin T \atop x\sim u}l(x)+ l(v)$,

$ Y_{v}\neq Y_{u} $.
\end{center}

Thus, $ Y_{v_{1}},\ldots ,Y_{v_{\omega}} $ are distinct numbers. On the other hand, for each vertex $ v \in T $, the image of the function $Y_{v}$ is $ [ -k , k(n -w )-1] $. So $ w\leq k(n -w+1) $, therefore $ k \geq \frac{w}{n-w+1} $ and the proof is completed.

{\bf (ii)}
Let $G$ be a regular graph, obviously $ \eta (G) \geq  2 $.
To the contrary suppose that $  \eta (G) = 2 $. Let $ T $ be a maximum clique in $ G $ and
$ c:V(G)\rightarrow \lbrace 1,2 \rbrace $ be an additive labeling of $ G $. Define:

\begin{center}
$ X_{1}=c^{-1}(1) \cap T $,\hspace{1cm}  $X_{2}=c^{-1}(2) \cap T  $,

$ Y_{1}=c^{-1}(1) \setminus T $,\hspace{1cm}  $Y_{2}=c^{-1}(2) \setminus T  $.
\end{center}

Suppose that $ X_{1}=\lbrace v_{1},\ldots ,v_{k}\rbrace $ and $ X_{2}=\lbrace v_{k+1},\ldots ,v_{\omega}\rbrace $.
For each $ i \in \mathbb{N}_{\omega} $, denote the number of neighbors of $ v_{i} $, in $ Y_{1} $ by $ d_{i} $.
Since $ c  $ is an additive labeling of the regular graph, every two adjacent vertices have different numbers of neighbors in $ c^{-1}(1) $. Therefore  $ d_{1},\ldots, d_{k},1+d_{k+1},\ldots ,1+d_{\omega}  $ are distinct numbers. Since for each $ i \in \mathbb{N}_{\omega}$, $ 0 \leq d_{i} \leq \vert Y_{1}\vert $, we have $\vert Y_{1}\vert \geq   \omega -2 $. Similarly, $\vert Y_{2}\vert \geq   \omega -2 $, so

\begin{center}
$ n=\vert T  \vert + \vert Y_{1} \vert + \vert Y_{2} \vert \geq 3 \omega -4$.
\end{center}

This is a contradiction. So the proof is completed.

}\end{alii}

\section{List Coloring Problem}
 \label{T4}

\begin{aliii}{

Let $G$ be a graph and let $L$ be a function which assigns to each vertex $v$ of $G$ a
set $L(v)$ of positive integers, called the list of $v$. A proper vertex coloring $c:V(G)\rightarrow \mathbb{N}$ such that $f(v)\in L(v)$ for all $v \in V $ is called a {\it list coloring} of $G$ with respect to $L$, or an {\it $L$-coloring}, and we say that $G$ is {\it $L$-colorable}.

Next, for a given graph $G$ and a list $L(v)$ for every vertex $v$, we construct a graph $H_G$ such that $H_G$ has a $(0,1)$-additive labeling if and only if $G$ is $L$-colorable.

Define $W= \bigcup_{v\in V(G)} L(v)$ and let $f$ be a bijective function   from the set $W$ to the set $\mathbb{N}_{|W|+1}\setminus \{1\}$. For every vertex $v\in V(G)$, let $L_{f}(v)=\{ f(i) | i \in L(v)\} $. The graph $G$ is $L$-colorable if and only if $G$ is $L_f$-colorable. Now, we construct $H_G$ form $G$ and $L_f$.

\textbf{Construction of  $H_G$}.\\
We use three auxiliary graphs $T(w)$, $I(j)$ and $G(v,L_f(v),s)$. The gadgets $I(j)$ and $T(w)$  are shown in Figure 1. Consider a vertex $v$ and a copy of auxiliary graph $T(w)$. Join the vertex $v$ to $T(w)$. Next, for every $j\in (\mathbb{N}_s\setminus \{1\}) \setminus L_f(v)$
consider a copy of $I(j)$ and join the vertex  $v$ to the vertex  $u_j$. Finally, put $s$ isolated vertices and join each of them to the vertex  $v$. Call the resulting graph $G(v,L_{f}(v),s)$. Now, for every vertex $v\in V(G)$ put a copy of $G(v,L_{f}(v),| W | +1)$ and for every edge
$v v' $ in the graph $G$ join the vertex $v \in V(G(v,L_{f}(v),| W | +1))$ to the vertex  $v' \in V(G(v',L_{f}(v'),| W | +1))$. Call the resulting graph $H_G$.

For a family $\mathscr{F}$ of graphs, define: $ \mathscr{F'} \stackrel {\rm def}{=}\{ H_G | G\in \mathscr{F}  \}$.
We  show that if
$\mathscr{F}$ is a family of graphs such that {\it list coloring  problem} is $ \mathbf{NP} $-complete for that family. Then, the following problem is $ \mathbf{NP} $-complete: "Given a graph $H_G \in \mathscr{F'}$, does $H_G$ have a $(0,1)$-additive labeling?"

\begin{figure}[h]
\begin{center}
\includegraphics[scale=.5]{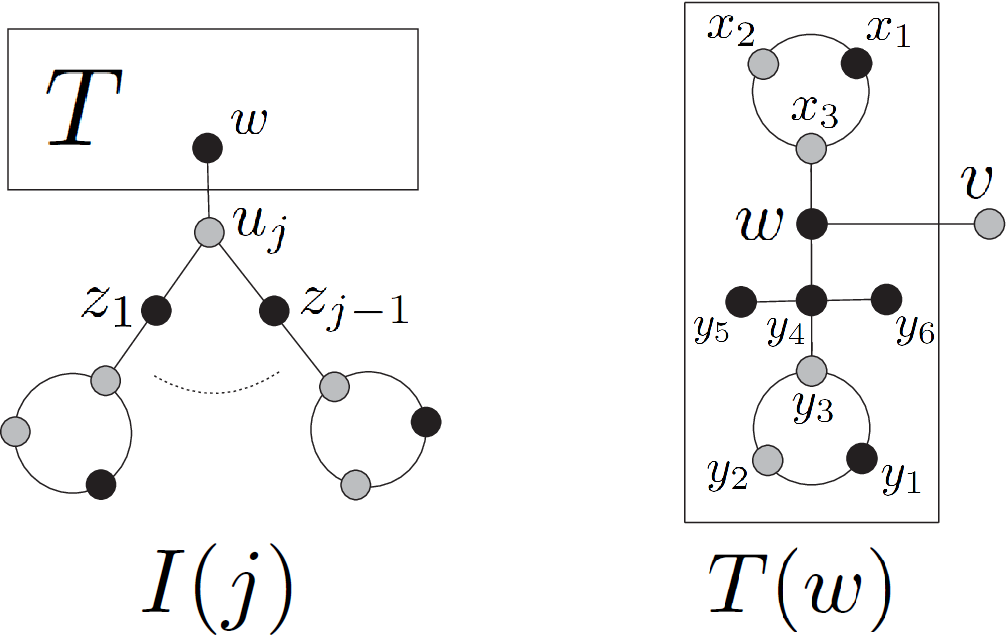}
\caption{The auxiliary graphs  $I(j)$ and $T(w)$.}
\end{center}
\end{figure}

First consider the following facts.

\begin{example}
\label{o1}
Let $G$ be a graph with a $(0,1)$-additive labeling $\ell$   and assume that it has the auxiliary graph $T(w)$ as a subgraph, $\ell(v) =0$, $\ell(w)=1$ and $\sum_{x\in N(w)}\ell(x) =1$.
\end{example}

{{\bf Proof of Fact \ref{o1}.}} By attention to the two triangles $x_1x_2x_3$ and $y_1y_2y_3 $, $\ell(w)=1$ and $\ell(y_4)=1$. Also $\ell(x_1)\neq \ell(x_2)$, without loss of
generality suppose that $\ell(x_1)=1 $ and $\ell(x_2)=0$. Therefore, $\ell(x_3)=0$, thus $\sum_{x\in N(w)}\ell(x) =1 + \ell(v)$. Since $ \sum_{x\in N(x_3)}\ell(x)=2$, therefore $\sum_{x\in N(w)}\ell(x) =1$, consequently
$\ell(v)=0$. $ \spadesuit$

\begin{example}
\label{o2}
Let $G$ be a graph with a $(0,1)$-additive labeling $\ell$   and assume that it has the auxiliary graph $I(j)$ as a subgraph,
 $\sum_{x\in N(u_j)}\ell(x) \geq j$.
\end{example}

{{\bf Proof of Fact \ref{o2}.}} By Fact \ref{o1}, $\ell(w)=1$, by using a similar argument  $\ell(z_1)= \cdots =\ell( z_{j-1})=1$. So
$\sum_{x\in N(u_j)}\ell(x)   \geq j$. $ \spadesuit$

\begin{example}
\label{F2}
Let $\ell$ be a $(0,1)$-additive labeling for $G(v,L_f(v),| W | +1)$, $\sum_{x\in N(v)}\ell(x)
 \in  L_f(v)$.
\end{example}

{{\bf Proof of Fact \ref{F2}.}}  By Fact \ref{o1} and Fact \ref{o2} it is clear.

First, suppose that the graph $H_G$ has a $(0,1)$-additive labeling $\ell$, define $c: V(G)\rightarrow \mathbb{N}$, $c(v)=\sum_{x\in N(v)}\ell(x)$. The function $c$ is a proper vertex coloring and for every vertex $v$, by Fact \ref{F2}, $c(v)\in L_f(v)$. Next, suppose that the graph $G$ is $L_f$-colorable, then it is clear that the graph $H_G$ has a $(0,1)$-additive labeling.

The list coloring problem is $ \mathbf{NP} $-complete for perfect graphs and planar graphs (see \cite{AAc}).
Obviously if $G$ is a planar graph, then $H_G$ is a planar graph.
Also, if $G$ is a perfect graph, then it is easy to see that the graph $H_G$ is a perfect graph.
This completes the proof.

}\end{aliii}

\section{Planar graphs}
 \label{T2}

\begin{aliiii}{

Let $\Phi$ be a $3$-SAT formula with the set of clauses $C $ and the set of variables
$X  $. Let $G(\Phi)$ be a graph with the vertices $C \cup X \cup (\neg X)$, where $\neg X = \lbrace \neg x: x\in X\rbrace$, such that for each clause $c=y \vee z \vee w $, $c$ is adjacent to $y,z$ and $w$, also every $x \in X$ is adjacent to $\neg x$. $\Phi$ is called planar $3$-SAT(type $2$) formula if $G(\Phi)$ is a planar graph. It was shown that the problem of satisfiability of planar $3$-SAT(type $2$) is $ \mathbf{NP}$-complete \cite{zhu1}. In order to prove our theorem, we reduce the following problem to our problem.

\textbf{Problem}: {\em Planar $3$-SAT(type $2$).}\\
\textsc{Input}: A planar $3$-SAT(type $2$) formula  $ \Phi $.\\
\textsc{Question}: Is there a truth assignment for $ \Phi $ that satisfies all the clauses?\\

Consider an instance of planar $3$-SAT(type $2$) with the set of variables
$X$ and the set of  clauses $C$. We transform this into a graph $G'(\Phi)$ such
that $ G'(\Phi) $ has a $(0,1)$-additive labeling, if and only if $\Phi$ is satisfiable.
The graph $G'(\Phi)$ has a copy of
$ B(x) $ for each variable $x$ and a copy of $A(c)$ for each clause $c$. The gadgets $ B(x)$ and $ A(c) $ are shown in Figure 2. Also, for every $c\in C$, $x\in X $, the
edge $w_{c}^{1}x$ is added if $c $ contains the literal $x$.
Furthermore, for every $c\in C$, $\neg x\in  \neg X $, the
edge $w_{c}^{1}\neg x$ is added if $c $ contains the literal $\neg x$.
Call the resulting graph $G'(\Phi)$.
Clearly the graph $G'(\Phi)$ is triangle-free and planar.

\begin{figure}[h]
\begin{center}
\includegraphics[scale=.4]{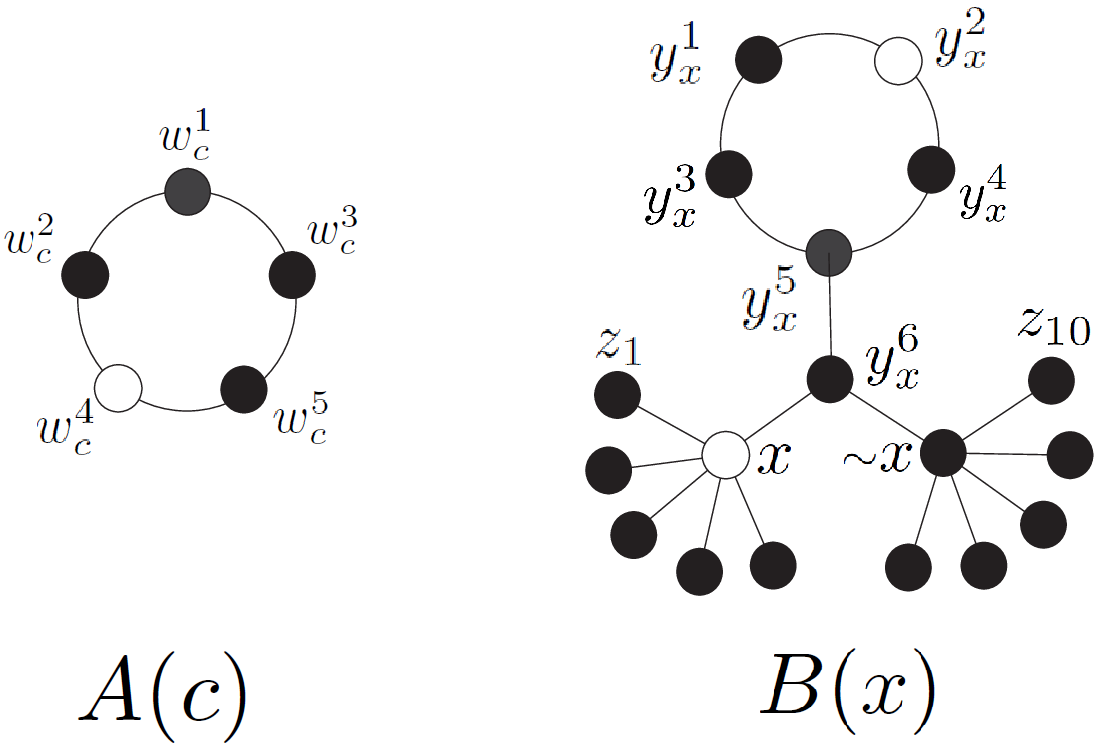}
\caption{The two auxiliary graphs $ A(c) $ and $ B(x)$.}
\end{center}
\end{figure}

\begin{example}
\label{F3}
Let $\ell$ be a $(0,1)$-additive labeling for the graph $G'(\Phi)$, for each clause $c=a \vee b\vee d$, $\ell(a)+\ell(b)+\ell(d)\geq 1$.
\end{example}

{{\bf Proof of Fact \ref{F3}.}} To the contrary suppose that there exists a clause $c =a \vee b\vee d$, such that $\ell(a)+\ell(b)+\ell(d)=0$, then $ \sum_{t\in N(w_{c}^{1})}\ell(t) = \ell( w_{c}^{2})+\ell(w_{c}^{3}) $.
Notice that in that case, $\ell$ restricted to the
odd cycle
$w_{c}^{1}w_{c}^{2}w_{c}^{4}w_{c}^{5}w_{c}^{3}$,
is a (0,1)-additive labeling, but an odd cycle
 does not have any $(0,1)$-additive labeling, this is a contradiction. $\spadesuit$

\begin{example}
\label{F4}
Let $G'(\Phi)$ be a graph with a $(0,1)$-additive labeling $\ell$, for each variable $x $, $\ell(x )+\ell(\neg x )\leq 1$.
\end{example}

{{\bf Proof of Fact \ref{F4}.}} To the contrary, suppose that there is a variable $x $, such that $\ell(x )+\ell(\neg x )=2$. Consider the auxiliary graph $B(x )$. Because of the odd cycle $y_{x}^{1} y_{x}^{2} y_{x}^{4} y_{x}^{5} y_{x}^{3}$, $\ell(y_{x}^{6})=1$. Now two cases for $\ell(y_{x}^{5})$ can be considered.

{{\bf Case $1$}}. $\ell(y_{x}^{5})=1$. Thus $\sum_{t\in N(y_{x}^{6})}\ell(t)  =3$, therefore $\sum_{t\in N(y_{x}^{5})}\ell(t) \in \{ 1,2\}$.

\textbullet $ $ If $ \sum_{t\in N(y_{x}^{5})}\ell(t)= 1$, then $\ell(y_{x}^{3})=\ell( y_{x}^{4})=0$. Thus,  $\ell(y_{x}^{1})+ \ell(y_{x}^{2})=1$; without loss of generality suppose that $\ell(y_{x}^{1})=1$ and $\ell(y_{x}^{2})=0$, in this case $ \sum_{t\in N(y_{x}^{2})}\ell(t)=\sum_{t\in N(y_{x}^{4})}\ell(t)$, but this is a contradiction.

\textbullet $ $ If $ \sum_{t\in N(y_{x}^{5})}\ell(t)= 2$. Suppose that $\ell(y_{x}^{3})=1$, $\ell(y_{x}^{4})=0$. Four subcases for $\ell(y_{x}^{1}),\ell(y_{x}^{2}) $ can be considered, each of them produces a contradiction.

{{\bf Case $2$}}. $\ell(y_{x}^{5})=0$. Thus $ \sum_{t\in N(y_{x}^{6})}\ell(t)= 2$, therefore $\sum_{t\in N(y_{x}^{5})}\ell(t) \in \{ 1,3\}$.

 \textbullet $ $ If $\sum_{t\in N(y_{x}^{5})}\ell(t) = 1$, then $\ell(y_{x}^{3})=\ell( y_{x}^{4})=0$. Therefore,  $\ell(y_{x}^{1})+\ell( y_{x}^{2})=1$. With no loss of generality suppose that $\ell(y_{x}^{1})=1$, $\ell(y_{x}^{2})=0$, therefore $ \sum_{t\in N(y_{x}^{3})}\ell(t)=\sum_{t\in N(y_{x}^{5})}\ell(t)$, but this is a contradiction.

 \textbullet $ $ If $ \sum_{t\in N(y_{x}^{5})}\ell(t) = 3$, then $\ell(y_{x}^{3})+\ell(y_{x}^{4})=2$. Thus $\ell(y_{x}^{1})+\ell( y_{x}^{2})=1$. Suppose that $\ell(y_{x}^{1})=1$, $\ell(y_{x}^{2})=0$, therefore $\sum_{t\in N(y_{x}^{1})}\ell(t)= \sum_{t\in N(y_{x}^{3})}\ell(t)$, this is a contradiction. $\spadesuit$

First, suppose that $ \Phi $ is satisfiable with the satisfying
assignment $ \Gamma: X \rightarrow \{ true, false\} $. We present a $(0,1)$-additive labeling $\ell$ for $G'(\Phi)$. For every variable $x $ if $\Gamma(x )=true $, then put $\ell(x )=1 $, otherwise put $ \ell(\neg x  )=1 $. Also put $\ell(  z_1)= \cdots= \ell(z_{10})=\ell(y_{x}^{1})=\ell( y_{x}^{3})=\ell( y_{x}^{4})=\ell( y_{x}^{5})=\ell( y_{x}^{6})=1$. Moreover, for every clause $c $, put $\ell(w_{c}^{1})=\ell( w_{c}^{2})=\ell( w_{c}^{3})=\ell( w_{c}^{5})=1 $. It is easy to extend this labeling to a $(0,1)$-additive labeling for the graph $G'(\Phi)$.
Next, suppose that the graph $G'(\Phi)$ has a $(0,1)$-additive labeling $\ell$. For each variable $x $, by Fact \ref{F4}, $\ell(x )+\ell( \neg x )\leq 1$. If
$\ell( x  )=1  $, put $\Gamma( x )=true $, if  $\ell( \neg x )=1
$, then put $\Gamma( x )=false $ and otherwise put $\Gamma( x )=true $. By Fact \ref{F3}, $\Gamma$ is a satisfying assignment for $ \Phi $.
}\end{aliiii}

\section{Inapproximability}
 \label{T3}

 \begin{aliiiii}{

Let $\varepsilon > 0$ and $ k$ be a sufficiently large number. It was shown that $3$-colorability of
$4$-regular planar graphs is $ \mathbf{NP} $-complete \cite{MR573644}. We
reduce this problem to our problem. In other words, for a given
$4$-regular planar graph $G$ with $k$ vertices, we construct a
planar graph $ G^*$ with $ 7k+10 k^{\lceil \frac{3}{\varepsilon}
\rceil +2} $ vertices, such that if $ \chi(G) \leq 3 $, then
$  \sigma_{1} (   G^*)\leq 5k$, otherwise $
\sigma_{1} (  G^*) >   5 k^{\lceil \frac{3}{\varepsilon} \rceil+1}$, therefore there is no $ \theta $-approximation algorithm
for determining $ \sigma_{1} (  G^*) $ for planar graphs, where:

\begin{center}
\begin{eqnarray*}
\theta =\frac{Approximate\ Answer}{OPT} & > & \frac{ 5 k^{\lceil \frac{3}{\varepsilon} \rceil+1}}{5k}\\
 & = &  k^{\lceil \frac{3}{\varepsilon} \rceil}\\
   & = &  \big( k^{\lceil \frac{3}{\varepsilon} \rceil +3} \big) ^{\frac{\lceil \frac{3}{\varepsilon} \rceil }{\lceil \frac{3}{\varepsilon} \rceil +3}}\\
 & \geq &   \big(7k+ 10k^{\lceil \frac{3}{\varepsilon} \rceil +2} \big) ^{\frac{\lceil \frac{3}{\varepsilon} \rceil }{\lceil \frac{3}{\varepsilon} \rceil +3}} \\
 & \geq &  \vert V(G^*)\vert ^{\frac{\lceil
\frac{3}{\varepsilon} \rceil}{\lceil \frac{3}{\varepsilon} \rceil +3}} \\
& \geq & \vert V(G^*)\vert ^{1-\varepsilon}
\end{eqnarray*}
\end{center}

In order to construct the graph $G^*$, we use the   auxiliary graph    $D(v)$  which is shown in Figure 3. Using simple local replacements, for every vertex $v$ of the graph
$G$, put a copy of $D(v)$, and for every edge $v u$ of the graph $G$,
join the vertex $v$ of $D(v)$ to the vertex $u$ of $D(u)$. Call the resulting graph $G^*$.
First, suppose that $G$ is not $3$-colorable and let $\ell$ be a $(0,1)$-additive labeling for $G^*$.
By the structure of $D(v)$   we have $\ell(v)=1$ and $\ell(p_3)=0$, so
$\sum_{x\in N(v)}\ell(x) = 4 + \ell(p_4)+\ell(p_5 )+\ell(p_6) $.
Since $G$ is not $3$-colorable, there exists a vertex $v$ such that $\sum_{x\in N(v)}\ell(x)=4$, therefore in the subgraph $D(v)$, $ \ell(p_4)+\ell(p_5 )+\ell(p_6) =0$, so $\ell(p_5)=0  $. Consequently for every $i$, $1 \leq i \leq d$, in the subgraph $D(v)$,  $\ell(v_i)+\ell(v'_i)\geq 1$. So $
\sigma_{1} (  G^*)> 5 k^{\lceil \frac{3}{\varepsilon} \rceil+1}$.
Next, suppose that $ \chi(G) \leq 3 $. So $G$ has a proper vertex coloring $
c:V(G)\rightarrow \lbrace 1,2,3\rbrace $. For every vertex $v$
of $G$, if $ c(v)=1 $ put $\ell( p_4)=\ell( p_6)=0 $ and
$\ell(p_5)=1 $, else if $ c(v)=2 $ let $\ell( p_4)=0 $ and
$\ell( p_5)=\ell( p_6)=1$ and if $ c(v)=3 $ let $ \ell( p_4)=\ell( p_5)=\ell( p_6)=1$. It is easy
to extend $\ell$ to a $(0,1)$-additive labeling for the graph $G^*$ such that $
 \sigma_{1} (  G^*)\leq 5k$.

\begin{figure}[h]
\begin{center}
\includegraphics[scale=.4]{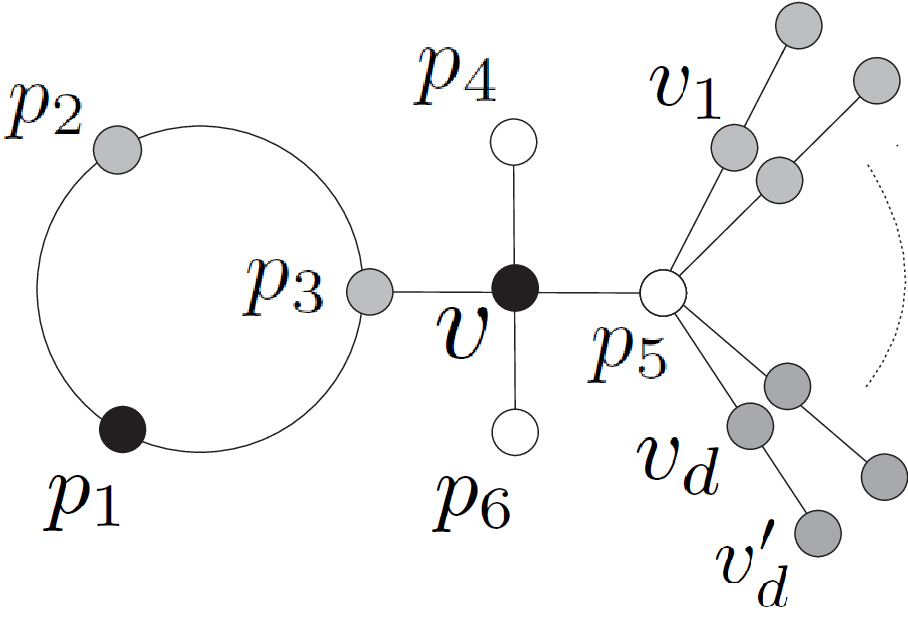}
\caption{The auxiliary graph $D(v)$. This graph has $7+ 10 k^{\lceil \frac{3}{\varepsilon}
\rceil +1}$ vertices, where $d=5 k^{\lceil \frac{3}{\varepsilon}\rceil +1}$. }
\end{center}
\end{figure}

}\end{aliiiii}

\section{Concluding remarks}
  \label{T5}

In this paper we study the computational complexity of $(0,1)$-additive labeling of graphs.
A
 $(0,1)$-additive labeling of a graph $G$ is a function $ \ell :V(G) \rightarrow\{0,1\}$, such that for every two adjacent vertices $ v $ and $
u$ of $ G $, $ \sum_{w \sim v}\ell(w)\neq \sum_{w \sim u}\ell(w) $. For future work, someone can consider another version of this problem that we call proper total dominating set. {\it A proper total dominating set} of a graph $G=(V,E)$, is a subset $D$ of $V$ such that every vertex has a neighbor in $D$ (all vertices in the graph including the vertices in the dominating set have at least one
neighbor in the dominating set) and every two adjacent vertices have a different number of neighbors in $D$ (note that in
a (0,1)-additive labeling every vertex does not need to have a neighbor labeled 1).

In this work, we proved that for every $k$ there is a graph $G$ such that $  \eta(G) \leq k \leq \eta_{\ell}(G)/2$.
What can we say about the difference in bipartite graphs?

\section{Acknowledgment}
\label{}

The authors wish to thank Saieed Akbari who drew their attention to Additive List Coloring Conjecture.

\bibliographystyle{plain}
\bibliography{luckyref}

\end{document}